\documentstyle[12pt]{article} 
\begin{document}
\setcounter{page}{1} 
\begin{center}
{\Large\bf Convolution numbers of binary sequences}
\medbreak\noindent
\centerline{by}
\medbreak\noindent
\centerline{Gregory M Constantine}
\centerline{Department of Mathematics}
\centerline{University of Pittsburgh}
\centerline{Pittsburgh, PA 15260}
\vskip.3cm\noindent
\centerline{and}
\vskip.3cm\noindent
\centerline{Rodica R A Constantine}
\centerline{Department of Psychology}
\centerline{University of Nevada}
\centerline{Las Vegas, NV 89119}
\vskip.7cm\noindent
\end{center}
\centerline{\bf ABSTRACT}
 
\medbreak\noindent Convolution sums are introduced and 
special instances of the cyclic convolution on finite 
sets is examined in more detail.  The distributions that 
emerge are multidimensional generalizations of the Catalan 
and Narayana numbers.  This work yields a closed form 
solution for 1-dimensional marginals and certain bivariate 
marginals in the cyclic prime case.  It is explained how 
a sufficiently high resolution of understanding these 
multidimensional distributions yields an approach to 
attack the Hadamard matrix conjecture.  

\medbreak\noindent 
\medbreak\noindent {\em Key words and phrases:\/}  Group 
action, probability distribution, paths, autocorrelation, 
Hadamard designs  \medbreak\noindent {\em Proposed running head:\/}  Convolution numbers 
\vskip2cm 
\footnoterule\noindent Funded under NIH grants 
P50-GM-53789, RO1-HL-076157 and an IBM shared University 
Research Award

\noindent gmc@pitt.edu 
  
\newpage
\section*{\Large\bf1.  Convolution sums} \medbreak\ Let $G$ 
be a group acting on a finite set $X$, and let 
$R$ be a ring.  We let $S$ be a set of functions from $X$ to 
$R.$ For $\sigma\in G$ and $f,g\in S$ define the element in $R$ 

\[\sigma (f,g):=\sum_{x\in X}f(x)g(\sigma x),\]
and call it the $\sigma -$convolution of $f$ with $g$. List the elements of $G$ in some arbitrary but fixed order; when writing $\sigma\in G$ we assume that $\sigma$ varies over the elements of $G$ in this specific order. Consider the vector
$G(f,g)=(\sigma (f,g):\sigma\in G).$ Partition $S\times S$ by $(f_
1,g_1)\sim (f_2,g_2)$ 
if $G(f_1,g_1)=G(f_2,g_2).$ Denote by $\{S_i\}$ the resulting 
equivalence classes.  Further write $G_i=G(f_i,g_i),$ where 
$(f_i,g_i)$ is a class representative of $S_i.$ Let $P$ be a 
probability measure on $S\times S.$ Induce a probability 
measure on $\{G_i\}$ by assigning $P(G_i)=P(S_i).$ The numbers 
$P(G_i)$ are called {\em convolution numbers.}
\medbreak\noindent A case of 
particular interest occurs when we restrict attention to 
just $\sigma (f,f),$ which we simply write as $\sigma (f);$ this yields 
an autocorrelation.  We abbreviate $G(f,f)$ by $Gf.$ In this 
case $P$ is a measure on $S$ and the $S_i$ are a partition of $S.$ 
The overarching goal is to understand this 
induced probability on $\{G_i\}$, for various $G,X,S$ and $R.$ We refer
to [7] for a thorough study of autocorrelation and its many applications
to engineering, probability and other applied branches of science. Studies of
autocorrelation are implicit in proofs of existence of codes and combinatorial
designs as is found in [11] and [9].

\medbreak\section*{\Large\bf2.  Definitions and notation in 
the cyclic case} \medbreak Of initial interest is the case 
of $G=C_n,$ the cyclic group of order $n$, and $X=C_n$ with 
cardinality as measure.  We take $R$ to be the ring $Z$ of 
integers.  Let $S$ be the set of functions from $X$ to $Z$ 
that take value 1 on exactly $k$ elements of $X$ and value 0 
on exactly $n-k$ elements of $X.$ Elements of $S$ are called 
0,1-{\em binary functions of weight} $k$ {\em on} $X$.  Clearly such 
functions are in bijection with subsets of $X$ with $k$ 
elements; evidently the cardinality of $S$ is $|S|={n\choose k}.$ We 
let $G=C_n$ act on $X=C_n$ by counterclockwise rotations, as is specified 
below.   
It is contextually clear when we interpret $C_n$ as a group 
or as the set on which the group acts, and this 
distinction is not always explicitly highlighted.  We also 
restrict to the autocorrelated case.  [In most, but not 
all, situations we may assume, without loss, that $k\leq\frac n2.
]$ 
\medbreak\noindent Write 
$G=C_n=\{\sigma_0=0,\sigma_1=1,\ldots ,\sigma_{n-1}=(n-1)\}.$ Then $
\sigma_i$ may be 
interpreted as a rotation of 'distance' $i$, in the sense that 
$\sigma_ix=x+i$ (mod $n$), for all $x\in X=C_n$.  The (periodic) 
autocorrelation $\sigma (f)=\sum_{x\in X}f(x)f(\sigma x)$ may then be viewed 
as the number of incidences of $f$ at 'distance' $\sigma$.  Put a 
uniform measure on $S.$ Write 
$Gf=(\sigma_0(f),\ldots ,\sigma_{n-1}(f)).$ Clearly $\sigma_0(f)=
k,$ for all $f;$ and 
$\sigma_i(f)=\sigma_{n-i}(f),$ for $i\geq 1.$ It thus suffices to write 
$Gf=(\sigma_1(f),\ldots ,\sigma_m(f))$, where $m=\lfloor\frac n2\rfloor$ and $
\lfloor x\rfloor$ stands for 
the integer less than or equal to the rational number $x$.  
\smallskip\noindent By writing $c(n,k;(d_1,\ldots ,d_m))$ for the 
the number of $k-$subsets of $X$ with $d_i$ incidences at 
distance $i,$ we may express 
\[P(Gf)={n\choose k}^{-1}c(n,k;(\sigma_1(f),\ldots ,\sigma_m(f)))
.\]
A couple of observations:  
\smallskip\noindent

{\bf 1.}  The number $c(n,k;(d_1,\ldots ,d_m))=0,$ unless 
$(d_1,\ldots d_m)=(\sigma_1(f),\ldots ,\sigma_m(f))$ for some $f\in 
S$.  Exactly 
when this number is 0 is a central issue in our 
undertakings, as is explained in some detail in Section 6.  
In essence, it is typically possible to find three binary 
vectors that could be included in the Goethals-Seidel 
construction.  These three vectors uniquely determine 
what the autocorrelation of the forth binary vector 
ought to be.  The Hadamard matrix construction can only 
be completed, however, if there exists a binary vector 
with this latter autocorrelation.  \vskip.3cm\noindent

{\bf 2.}  We have $\sum_{i=1}^md_i={k\choose 2},$ since for any $
k-$subset we 
have ${k\choose 2}$ unordered distances available in all.  
\vskip.3cm\noindent An upward-and-right moving path on 
the integer lattice starting at (0,0) and ending at $(n-k,k)$ 
that touches or stays above the line joining (0,0) and 
($n-k,k)$ is, for simplicity, just called a $path.$ Note first 
that any path may be viewed as a $k-$subset (and, 
equivalently, as a 0,1-binary sequence with $k$ ones), by 
listing the indices of the upward moves as elements of 
the set. To be specific, an upward move is marked by 1 and a move to the right by 0. We draw further attention to the following 
useful observation.  \medbreak\noindent {\bf Lemma} {\em The $C_n$-orbit
of any binary sequence of length $n$ and weight $k$ contains at least 
one path. It contains exactly one path if $n$ and $k$ are coprime.} 
\medbreak\noindent
{\bf Proof} This is best seen as follows. The line joining $(0,0)$ 
and $(n-k,k)$ has slope $\frac {k}{n-k}$. For any 0-1 sequence 
$s=(s_i:1\leq i\leq n)$ of length $n,$ replace 
the 0s by $\frac {-k}{n-k}$ and leave the 1s as they are; call the new 
sequence $s'=(s_i').$ Calculate the partial sums $p_i$ of $s'$ by 
setting $p_1=s_1',$ and $p_i=p_{i-1}+s_i'$, $2\leq i\leq n.$ Let $
i^{*}$ be an 
index for which $p_{i^{*}}$ is a minimum. Apply a cyclic 
rotation to the original sequence $s$ that places in 
position 1 the index $i^{*}+1.$ It is evident from this 
construction, by the choice of $i^{*}$, that the resulting 
sequence is in the $C_n$-orbit of $s$ and that it corresponds 
to a path. When $k$ and $n$ are coprime the index $i^{*}$ is 
unique; else the path would touch the diagonal line at 
$(h,v)$ with $(h,v)\neq (0,0)$ or $(n-k,k)$. Similar right triangles 
now yield $\frac vh=\frac k{n-k},$ or $vn=k(v+h).$ Since $k$ and $
n$ are 
coprime this forces $k$ to divide $v,$ but this is not 
possible because $v<k.$  This ends the proof.  
\vskip.3cm\noindent
[For instance, if $s=(0,1,0,0,1,1,0),$ then its orbit is 
represented by the shifted sequence (1, 1, 0, 0, 1, 0, 0) which 
corresponds to a path.  Note that $s$ itself does not 
correspond to a path.]  \medbreak\noindent By a {\em descent\/} in 
a 0,1-binary sequence we mean an occurance of 10 in the 
sequence.  When the binary sequence corresponds to a 
path, a descent occurs when an upward move is followed 
by a right move on the path.  [It is easy to see, for 
example, that the sequence (0, 1, 0, 0, 1, 1, 0) has exactly two 
descents.] An $ascent$ is an occurance of 01 in the sequence.  
\vskip.3cm\noindent
A {\em block of ones\/} in a 0,1-binary sequence is a 
subsequence of consecutive 1s; it is called a {\em run of ones}
if it is maximal, by inclusion, with this property.  A 
{\em run of zeros} is analogously defined.  Observe that a run 
of ones always ends in a descent and a run of zeros 
always ends in an ascent.  

\medbreak\section*{\Large\bf3.  The one-dimensional 
marginals} \medbreak\medbreak Define $D=(D_1,\ldots ,D_m)$ to 
be a random vector, with $D_i$ taking values in the set 
$\{0,1,\ldots ,k\}.$ The component $D_i$ counts the number of 
incidences at distance $i$ that can occur for an arbitrary 
$k-$subset.  Describing explicitly the joint probability 
$P(D=(d_1,\ldots ,d_m))$ seems difficult.  Tractable are the one 
and possibly two-dimensional marginal distributions of $D.$ 
We abide by the usual notational conventions when 
working with binomial numbers but also draw attention 
to the fact that we convene to write ${t\choose t}=1$ for all 
integers $t$; thus ${{-3}\choose {-3}}=1.$ 

\medbreak\noindent {\bf Theorem}
\vskip.3cm\noindent
{\em $(a)$ If} $i$ {\em and} $n$ {\em are coprime, then} $D_i$ {\em has the following}
{\em distribution that does not depend on} $i${\em :}
\[{n\choose k}P(D_i=x)=\frac n{(k-x)}{{k-1}\choose {k-x-1}}{{n-k-
1}\choose {k-x-1}},\]
{\em for} $x=0,1,\ldots ,k-1${\em ;} $k\leq\frac n2,n\geq 3.$

\smallskip {\em $(b)$ Let} $a(v,s;j,x)$ {\em denote the number of}
$s-${\em subsets of a set with} $v$ {\em elements that have} $x$ 
{\em incidences at distance} $j;$ {\em define} $a(v,0;1,0)=0.$ {\em If} $
i(>1)$ 
{\em divides} $n,$ {\em we obtain the distribution of} $D_i$ {\em inductively}
{\em by writing}
\[a(n,k;i,x)={n\choose k}P(D_i=x)=\sum_{{\bf w}}(\sum_{k_j\in {\bf w}}
a(\frac ni,k_j;1,x)),\]
{\em with the sum ranging over all vectors} ${\bf w}=(k_1,\ldots 
,k_i)$ {\em with}
$0\leq k_j\leq\frac ni,$ $\sum_jk_j=k$$.$ {\em Here} $x=0,1,\ldots 
,k.$ {\em Since} $i>1$, {\em the terms}
{\em in the sums on the right-hand-side are known by}
{\em induction and part $(a)$.}
\smallskip\

{\em $(c)$ The distribution of} $D_i$ {\em in part $(a)$ is log-concave, hence}
{\em unimodal. It increases in} $r=k-x$ {\em from 1 to} $\frac {k
(n-k)-1}n$ 
{\em and decreases afterwards;} $1\leq r\leq k.$
  
\medbreak\noindent
Note that for $n$ odd, $k=\frac {n-1}2$, and $i$ coprime to $n$, the 
above Theorem yields
\[{n\choose k}P(D_i=x)=\frac {2k+1}{k-x}{{k-1}\choose {k-x-1}}{k\choose {
k-x-1}}=nN(k,x),\]
with $N(k,x)$ being the Narayana numbers [1].  
\medbreak\noindent {\bf Proof} Since the automorphism group of 
$G=C_n$ acts transitively on the primitive roots (or 
generators) of $C_n$, it suffices to prove part $(a)$ for $i=1.$ 
Specifically, if $\phi$ sends $i$ to 1, then a $k-$subset $A$ with $
x$ 
incidences at distance $i$ is sent by $\phi$ into a $k-$subset 
$\phi (A)$ that has $x$ incidences at distance 1; this yields, 
$P(D_i=x)=P(D_{\phi (i)}=x)=P(D_1=x).$ \smallskip\noindent
Any cyclic 0,1-binary sequence with $k$ ones can always 
be listed so as to have the beginning of a run of ones in 
position 1.  Mark such a cyclic sequence as
\medbreak\noindent
\centerline{ $1^{a_1}0^{b_1}1^{a_2}0^{b_2}\cdots 1^{a_r}0^{b_r}=\prod_{
i=1}^r1^{a_i}0^{b_i}$,}
\medbreak\noindent
with $1^{a_i}0^{b_i}$ indicating that a run of $a_i$ ones is followed by a run of 
$b_i$ zeros; $a_i,b_i\geq 1.$  Evidently $\sum a_i=k$ and $\sum b_
i=n-k.$ We seek 
the number of $k-$subsets with $x$ incidences at distance 1.  
A run with $a_i$ ones yields $a_i-1$ such incidences and 
therefore $x=\sum (a_i-1)=(\sum a_i)-r=k-r,$ where $r$ denotes 
the number of runs of ones. Viewing the sequence in 
terms of runs of ones, it is 
apparent that the number we seek is equal to the number 
of interlaced compositions with $r$ parts (the run sizes $a_i$ 
of ones) of $k$, and the compositions with $r$ parts (the 
run sizes $b_i$ of zeros) of $n-k.$ This yields the product 
${{k-1}\choose {r-1}}{{n-k-1}\choose {r-1}}$ as an initial count. However, for a cyclic sequence we 
have $n$ places to choose from for an initial first index, so we 
multiply   
this product by $n$; we also must divide it by $r$ since we can only 
use one of the $r$ runs of ones to start in position 1. 
Using the fact that $x=k-r,$ part $(a)$ is now 
demonstrated. 
\vskip.3cm\noindent\ 
Part $(b)$ addresses the case when $i$ divides $n.$ We provide an 
inductive answer. View the $k-$subset as a binary 
sequence $s.$ Split $s=(s_j)$ into $i$ subsequences each of length $\frac 
ni$, 
the $m^{th}$ subsequence being $s_m,s_{m+i},s$$_{m+2i},\ldots ,s_{
m+(\frac ni-1)i};$ 
$1\leq m\leq i.$ The statement becomes clear upon observing 
that there are $k$ ones in $s$ if and only if the $k$ ones are 
partitioned in all ways possible among the $i$ 
subsequences. Incidences at distance $i$ in $s$ become 
incidences at distance 1 within each subsequence.
\vskip.3cm\noindent
Part $(c)$ can be verified through a direct calculation. If 
$p_r=\frac 1r{{k-1}\choose {r-1}}{{n-k-1}\choose {r-1}}$, $r=1,\ldots 
,k$ then the ratio 
\[\frac {p_{r+1}}{p_r}=\frac {(k-r)(n-k-r)}{r(r+1)}\]
shows that $p_r$ increases for $1\leq r\leq\frac {k(n-k)-1}n$ and decreases 
afterwards. Log concavity is also checked directly by 
verifying that $p_r^2\geq p_{r-1}p_{r+1}.$ This inequality is shown 
equivalent to $2k(n-k)\geq (r-1)(n+1)$, which suffices to be 
checked for $r=k,$ since $r\leq k$. Since $2k\leq n$, the inequality 
is true. 
\medbreak
\section*{\Large\bf4. Two-dimensional marginals} \medbreak
For $n$ and $k$ $(\leq\frac n2)$ we shall give an explicit formula for the 
number of $k-$subsets with $x$ incidences at distance 1 and 
$y$ incidences at distance 2. This is equivalent to finding 
$P(D_1=x,D_2=y).$ 
\newpage 

\medbreak\noindent {\bf Proposition} {\em The number of k-subsets 
that have x incidences at distance one and y incidences at distance two is equal to}
\[{n\choose k}P(D_1=x,D_2=y)=\]
\[\frac n{(k-x)}\sum_{a+b=k+y-2x}{{k-x}\choose a}{{x-1}\choose {k
-x-a-1}}{{k-x}\choose b}{{n-2k+x-1}\choose {k-x-b-1}}\]
{\bf Proof} As in the proof of the Theorem, summarize a 
$k-$subset in the notation $\prod_{i=1}^r1^{a_i}0^{b_i},$ which highlights the 
runs in the corresponding binary sequence.  In a run of 
size $a_i>1$ we have $a_i-2$ incidences at distance 2, if 
$b_i>1,$ and $a_i-1$ incidences at distance 2, if $b_i=1.$ It 
follows that 
\[y=\sum_{a_i>1}(a_i-2)+|\{i:b_i=1\}|=\]
\[\sum_i(a_i-2)+|\{i:a_i=1\}|+|\{i:b_i=1\}|=\]
\[k-2r+a+b=(k-r)-r+a+b=2x-k+a+b,\]
where $x=k-r,$ $a=$$|\{i:a_i=1\}|$ and $b=|\{i:b_i=1\}|.$
\smallskip\noindent
We note that $a+b$ counts the total numbers of runs of 
size 1 (be they of 1s -- which is a, or of 0s -- which is 
b).  By above, $a+b=k+y-2x.$ We conclude that {\em the }
{\em number of} $k-${\em subsets with} $x$ {\em incidences at distance 1 and}
$y$ {\em incidences at distance 2 is equal to} $\frac n{(k-x)}$ {\em times the}
{\em number of paths with} $r=k-x$ {\em one-runs (or descents)}
{\em and with} $a+b=k+y-2x$ {\em runs of size 1.\/}  We multiply 
by $n$ since paths are orbit representatives (as stated in 
the Lemma), and 
divide by 
$r=k-x$ since only one of the $r$ expressions for the 
$k-$subset in the form $\prod_{i=1}^r1^{a_i}0^{b_i}$ is the path in question.  
To count the number of such paths, express a path in 
the form $1^{a_1}0^{b_1}\cdots 1^{a_r}0^{b_r}$ in which $a$ of the $
a_i$s are 1 and 
$b$ of the $b_i$s are 1, with $a+b=k+y-2x.$ The vector 
$(a_1,\ldots ,a_r)$ is a composition of $k,$ and $(b_1,\ldots ,b_
r$) is a 
composition of $n-k.$ Focus on the composion of $k.$ Place a 
1 in each of its r parts; we now count unrestricted 
compositions of $k-a$ with $r-a$ classes; this yields 
${r\choose a}{{k-r-1}\choose {r-a-1}}$ as answer. Argue analogously for the 
composition on $n-k.$ This completes the proof.
\medbreak
\section*{\Large\bf5.  Examples and connection to Catalan 
numbers} \medbreak We begin with a motivational 
example. Take a set with $n=15$ elements and examine all 
its subsets of size $k=6.$ Our interest is in counting the 
number of (circular) incidences that occur at distance $i$ 
across all such subsets, as explained at the beginning of 
Section 2. Table 1 displays, by a direct count, the number 
of $k-$subsets with $x$ incidences at distance $i,$ for $i=1,3,5$ 
and $x=0,1,2,3,4,5,6.$ Provided that $j$ is coprime to 15, it 
is verified that the number of incidences at distance $j$ 
are equinumerous to those at distance 1; this is 
highlighted in part (a) of the Theorem. All entries in 
Table 1 are in agreement with the values provided by 
our Theorem. 
\begin{table}[ht]
\caption{}
\centering
\begin{tabular}{| c | c | c | c | c | c | c | c |}
\noalign{\vskip 2mm}
\hline
\hline
&\multicolumn{7}{|c|}{$x$}\\
\cline{2-8}
$i$ & 0 & 1 & 2 & 3 & 4 & 5 & 6 \\
\hline
1 & 140 & 1050 & 2100 & 1400 & 300 & 15 & 0 \\
\hline
3 & 125 & 1125 & 1950 & 1550 & 225 & 30 & 0 \\
\hline
5 & 0 & 1215 & 2430 & 810 & 540 & 0 & 10 \\
\hline
\hline
\end{tabular}
\label{table:nonlin}
\end{table}

\vskip.3cm\noindent
The numbers displayed in a row in Table 1 are 
${{15}\choose 6}P(D_i=x)$, as counted in our Theorem. It may be worth noticing that for 
$i=5$, a divisor of 15, the sequence of probabilities is not 
log-concave. We point out that the Theorem makes no statement on 
log-concavity in such cases.  \vskip.3cm\noindent We 
outline now in some detail how the entries in the row $i=5$ are obtained, 
using part (b) of the Theorem.  The split subsequences 
of indices are 1 6 11; 2 7 12; 3 8 13; 4 9 14; 5 10 15.  The 
possible compositions of $k=6$ (we just list partitions, 
for brevity) are 2 1 1 1 1, 3 1 1 1 0, 2 2 1 1 0, 3 2 1 0 0, 2 
2 2 0 0, and 3 3 0 0 0.  For instance, partition 2 2 1 1 0 
-- which carries two incidences -- indicates the fact 
that, of the $k=6$ available ones, we distribute 2 ones in 
two subsequences, 1 one in two subsequences, and 0 ones 
in one subsequence.  The number of ways of doing this 
is $5\cdot{3\choose 2}\cdot 4\cdot{3\choose 2}\cdot 3\cdot{3\choose 
1}\cdot 2\cdot{3\choose 1}/(2!\cdot 2!)=2430,$ which 
explains the entry corresponding to $x=2$ in the $i=5$ 
row.  The other entries are analogously explained; 
observe, for example, that $x=3$ arises from both 3 1 1 1 
0 and 2 2 2 0 0.  \medbreak We shall now explore a couple of 
consequences of the Theorem.  Assume that $n$ and $k$ are 
coprime and fixed, and let 
$a(i,x):=a(n,k;i,x)={n\choose k}P(D_i=x)$, with the latter having 
the explicit form written in the Theorem.  As $G=C_n$ 
acts on $k-$subsets, notice that coprimality implies that 
each orbit has length $n;$ this is either seen directly or 
can be a consequence of the Cauchy-Frobenius lemma, 
since no element except the identity has fixed points.  
This tells us that $\frac 1n{n\choose k}$ counts the number of $C_
n$-orbits 
and is therefore an integer.  It tells us also that $a(i,x)$ 
counts the lengths of certain types of orbits and is 
itself necessarily divisible by $n.$ By letting $b(i,x)=\frac {a(
i,x)}n$ 
and $b_{n,k}=\frac 1n{n\choose k}$, we obtain from the Theorem the 
following consequence.  \medbreak\noindent {\bf Corollary 1} {\em If 
$n$ and $k$ are coprime, then $b(i,x)$ counts the number of}
$C_n-orbits$ {\em that have $x$ incidences at distance $i$,}
$1\leq i\leq\frac {n-1}2.$ {\em These numbers satisfy} $\sum_xb(i
,x)=b_{n,k},$ {\em for all i, and provide a refinement of the numbers} $b_{n,k}.$ 
{\em When $i$ is coprime to $n$ the numbers $b(i,x)$ are all positive integers.  In particular, for n odd and} $k=\frac {
n-1}2${\em , }
{\em $b(1,x)$ are the Narayana numbers and} $b_{n,\frac {n-1}2}$ {\em is the Catalan number.}
\medbreak\noindent
For $n=21$ and $k=10$ 
we obtain the refinement of the Catalan number [2] as a 
sum of Narayana numbers as follows:  \hfill\break 
\hfill\break
\centerline{$1+45+540+2520+5292+5292+2520+540+45+1=16796$}
\[\]

In parallel, for $n=21$ and $k=8$ a refinement of $b_{21,8}$ in 
terms of the $b(1,x)$ given by Corollary 1 is\hfill\break

\centerline{$99+924+2772+3465+1925+462+42+1=9690$}

\medbreak\noindent {\bf Corollary 2} {\em If n is odd,} $k=\frac {
n-1}2${\em , and $i$}
{\em is coprime to $n$, then $b(i,x)$ are the Narayana numbers;}
{\em when $i$ is not coprime to $n$ the non-negative integers}
{\em $b(i,x)$ provide a refinement of the Catalan number} $b_{n,k}$ 
{\em that is different from the Narayana numbers.}
\medbreak\noindent For example, the Catalan number 
$b_{15,7}=429$ has the following three refinements, 
displayed in Table 2 below, offered by Corollary 2. The 
middle row, corresponding to $i=4,$ displays the 
Narayana sequence. 
\begin{table}[ht]
\caption{} 
\centering
\begin{tabular}{| c | c | c | c | c | c | c | c |}
\noalign{\vskip2mm}
\hline
\hline
&\multicolumn{7}{|c|}{$x$}\\
\cline{2-8}
$i$ & 0 & 1 & 2 & 3 & 4 & 5 & 6 \\
\hline
3 & 0 & 25 & 100 & 175 & 110 & 17 & 2 \\
\hline
4 & 1 & 21 & 105 & 175 & 105 & 21 & 1 \\
\hline
5 & 0 & 0 & 162 & 135 & 108 & 18 & 6 \\
\hline
\hline
\end{tabular}
\label{table:nonlin}
\end{table}
\medbreak\noindent\ 
In general, with the exception of the Narayana case, integer 
sequences $b(i,x)$ and $a(i,x)$ are not found in Sloane's 
Encyclopedia of Integer Sequences [6].
\medbreak

\section*{\Large\bf6.  A connection to the construction of 
Hadamard matrices} \medbreak A Hadamard matrix is a square matrix with entries -1 or 1 and orthogonal rows. It is easy to see that the dimension of such a matrix is either 2 or a multiple of 4. A remaining central question is whether these conditions are also sufficient for existence. A leading and remarkable result found in [8] informs us that if the dimension is divisible by a sufficiently high power of 2, then a Hadamard matrix of that dimension exists.
A detailed research resource on the subject appears in [9]. Sporadic yet useful constructions are found in [10]. Goethals and Seidel [4] give a 
method of using four circulant matrices, each of 
dimension $n$, to construct a Hadamard matrix of order 
$4n$ ([3] and [5]).  Each circulant can be specified by a binary (say 0,1) 
vector of length $n$; the sufficient condition which yields 
the Hadamard matrix is that the four binary 0,1 vectors 
$f_1,f_2,f_3,f_4$ (with $f_i$ of weight $k_i$) have (cyclic) autocorrelations 
$\sigma (f_i)$ that sum to the constant vector with all its entries equal 
to $k_1+k_2+k_3+k_4-n$.  
Vectors $\sigma (f_i)$ are of length $m=\frac {(n-1)}2$, as explained in 
Section 1.  \smallskip\noindent Any autocorrelation vector 
of a $k$-subset has nonnegative integral entries that sum 
to ${k\choose 2}$; but these necessary conditions are far from 
being sufficient.  The content of this paper can be 
viewed as an attempt to make headway toward 
establishing necessary and sufficient conditions in the 
form of understanding the joint distribution of the 
autocorrelation vector; and, in particular, in having a 
complete understanding of its (nonzero) support.  In 
other words, holding in our hand a vector with 
nonnegative entries that sum to ${k\choose 2}$ we should be 
able to tell with certainty whether or not this is the 
autocorrelation vector of some $k$-subset.  This paper 
established the one and two-dimensional marginal 
distributions of the autocorrelation.  Even this limited information 
proves (marginally) helpful, as we shall see in the examples 
below.  \smallskip\noindent A manageable case is n=19, 
with $(k_1,k_2,k_3,k_4)=(7,6,9,9)$.  We can take both $f_3$ and $
f_4$ 
to be 
the quadratic residue binary sequence that has 
autocorrelation vector of length 9 with all entries equal 
to 4.  We now seek a 7-subset $f_1$ and a 6-subset $f_2$ whose 
autocorrelation sequences sum to (7+6+9+9-19)-(4+4)=4 in each of the 9 
components. Such pairs of subsets are called 
{\em supplementary difference sets\/} and they are generally 
difficult to construct. 
\smallskip\noindent
Typically, starting with a 7-subset $f_1$ we attempt to 
supplement it with a 6-subset $f_2.$ There are obvious 
necessary conditions on $f_1$ but no complete understanding 
exits on how to complete this task. For instance, we 
may start with the 7-subset 
\[f_1=(0,0,1,0,1,0,0,1,0,1,0,0,1,0,0,1,0,1,0)\]
and attempt to supplement it. In this case 
\[\sigma (f_1)=(0,3,3,1,4,3,1,4,2);\]
This forces $\sigma (f_2)=(4,1,1,3,0,1,3,0,2),$ the supplement to 
the vector with all entries 4. But such a 6-subset $f_2$ 
does not exist. Indeed our Proposition tells us that there are 
no 6-subsets with 4 incidences at distance 1 and 1 
incidence at distance 2 (this can also be verified directly 
in this small case). 
\smallskip\noindent
But such a supplementary pair does exist. Were we to start 
with 
 
\[f_1=(0,0,0,1,1,0,0,1,0,1,0,0,1,0,0,1,1,0,0),\]
which has $\sigma (f_1)=(2,1,3,2,2,3,3,3,2),$ we would be seeking a supplement 
$f_2$ that must have $\sigma (f_2)=(2,3,1,2,2,1,1,1,2).$ Indeed, such $
f_2$ exists, an 
example being 
\[f_2=(1,1,1,0,1,0,1,0,0,0,0,1,0,0,0,0,0,0,0).\]
[As an aside that might make us aware of the delicate 
combinatorics involved: If we switch just entries 17 and 
18 in the vector $f_1$ above, a solution ceases to exist!]
\medbreak\noindent
We examine now a couple of larger cases. A Hadamard matrix of order 
4$\cdot$79 is known to exist. We use the Goethals-Seidel method to construct one here. Four binary vectors of length 79 and weights 34, 34, 42, 43 are needed, whose autocorrelations must sum to 74 in each of the 39 entries. Initially we build three such vectors of weights 34, 34, 42 the sum of autocorrelations of which differs from the vector with all entries equal to 74 by the vector 
\medbreak\noindent
$a=$ 22 23 25 25 26 24 25 21 22 23 24 22 24 21 24 21 24 26 25 24 24 21 25 25 22 24 21 22 21 26 20 21 23 23 23 21 25 23 22
\medbreak\noindent
Our Theorem and Proposition do not preclude the existence of a binary vector $v$ of length 79 and weight 43 that has $a$ as its autocorrelation. Indeed, a computer search yields such a vector.
The initial three vectors, coded in base 10 as 0 through 7, with column $(0,0,1)^{t}$ 
representing the number 1, are
\medbreak\noindent
4 3 3 2 7 1 0 1 7 4 6 3 1 0 0 0 1 5 7 3 6 1 1 1 6 5 5 4 3 4 7 4 6 2 0 2 7 0 5 4 0 1 0 3 2 0 6 6 4 3 0 3 4 1 1 2 5 1 5 2 7 7 1 5 4 0 0 1 3 2 0 7 5 0 5 7 2 7 3
\medbreak\noindent
The sought-after vector $v$, which allows us to complete the Goethals-Seidel construction in this case, is
\medbreak\noindent
$v=$ 1 1 1 0 1 1 0 1 1 0 0 1 1 0 0 0 1 1 1 1 0 1 1 1 0 1 1 0 1 0 1 0 0 0 0 0 1 0 1 0 0 1 0 1 0 0 0 0 0 1 0 1 0 1 1 0 1 1 1 0 1 1 1 1 0 0 0 1 1 0 0 1 1 0 1 1 0 1 1
\medbreak\noindent
The first open case for which a Hadamard matrix is not known to exist 
involves dimension 668. Using the same Goethals-Seidel method,
one way to attempt to construct one is to proceed as in the 79-case 
studied above and use four binary vectors, each of 
length 167, of weights 76, 76, 77, 80 whose four autocorrelations sum to 142 
in each of the 83=(167-1)/2 components. With some effort, it is possible to 
find three vectors of weights 76, 76, 77 whose sum of the three autocorrelations
differ from 142 in each component as shown in the vector $b$ written below. 
\medbreak\noindent
$b=$ 41 37 40 39 41 39 35 38 36 34 42 39 37 37 32 37 37 36 36 36 40 39 36 37 39 38 35 33 39 35 37 41 42 40 41 38 43
41 34 39 39 36 42 39 38 41 40 40 37 38 37 37 35 37 37 37 36 38 37 34 37 40 39 37 38 38 42 38 42 37 34 39 39 37
41 36 34 37 38 42 45 39 40
\medbreak\noindent
The three vectors, coded in base 10 as 0 through 7, as before, are as follows:
\medbreak\noindent
1 6 1 6 7 1 1 6 1 0 7 6 6 0 6 0 7 0 1 7 0 0 1 6 1 0 6 7 7 6 0 1 6 7 1 6 0 0 0 0 6 7 7 1 1 0 1 0 1 0 1 1 7 0 0
1 0 1 1 1 6 6 7 6 1 0 1 7 6 1 7 6 1 0 1 1 0 6 6 6 0 1 6 0 7 0 0 6 0 7 7 0 6 0 7 1 0 7 0 0 6 6 7 0 6 6 1 0 0 0
1 7 0 0 1 7 0 6 6 0 1 0 6 7 7 7 1 0 7 6 1 7 7 0 7 6 6 0 7 1 7 1 7 1 7 0 0 6 1 6 6 6 1 7 7 7 6 1 6 7 0 0 0 1 6
0 0
\medbreak\noindent
To have a successful construction one must verify that the vector $b$ is indeed 
the autocorrelation of a binary vector of length 167 and weight 80. Recognizing 
whether this is the case or not is the central motivational issue for this paper. 
Currently we simply do not know whether this is true or 
not, but this article does not center on this issue. Recent advances in algorithmic verification along 
these lines are found in [12].
\medbreak\noindent
The parting remark is that if the full joint distribution 
(or the support of positive probability) of the 
autocorrelation is understood, then we can decide whether or 
not a Hadamard matrix can be constructed by the 
Goethals-Seidel circulant method without necessarily explicitly providing the 
four required circulants. The approach may, in other words, assert existence without 
relying on explicit construction. 
\medbreak\noindent
\centerline{{\bf Acknowledgement}}
\medbreak\noindent
The article benefited from relevant clarifications made by one Referee, and for suggesting the inclusion of 
larger motivational examples of the Goethals-Seidel 
construction. 

\newpage

\centerline{{\bf REFERENCES}}
\begin{enumerate}

\item Narayana, T. V. Lattice Path Combinatorics with Statistical Applications. Toronto, Canada: University of Toronto Press, 100 - 101 (1979)

\item Stanley, R. P. {\em Enumerative Combinatorics}, Vol. 2. Cambridge, England: Cambridge University Press, 1999 

\item MacWilliams, F. J. and Sloane, N. J. A., {\em The theory of error-correcting codes}, Elsevier/North-Holland, 
Amsterdam, 1977

\item Goethals, J. M. and J. J. Seidel, Orthogonal matrices with
zero diagonal, Canad. J. Math., 19, 1001 - 1010 (1967)

\item de Launey, W. and D. Flannery, {\em Algebraic design theory}, Mathematical
Surveys and Monographs, vol. 175, American Mathematical Society,
Providence, RI, 2011

\item OEIS Foundation Inc. (2018), The On-Line Encyclopedia of Integer Sequences, http://oeis.org 

\item Damelin, S. and Miller, W. {\em The Mathematics of Signal Processing}, Cambridge University Press, 2011, ISBN 978-1107601048

\item Seberry, J., On the existence of Hadamard matrices, J Comb Theory (A), 21, 188-195 (1976)  

\item Seberry J., {\em Orthogonal designs}, Springer, Cham, 2017. Hadamard matrices, quadratic
forms and algebras, Revised and updated edition of the 1979 original [MR0534614]

\item Djokovic D., Golubitsky O., and Kotsireas I., Some new orders of Hadamard
and skew-Hadamard matrices. J. Combin. Des., 22(6):270-277 (2014)

\item Lander E., {\em Symmetric designs:  an algebraic 
approach}, Cambridge University Press, London, 1983 

\item Bright, C., Kotsireas, I. and Ganesh, V. Applying 
computer algebra systems with SAT solvers to the 
Williamson conjecture, {\em Journal of Symbolic Computation}, 
vol. 100, 187-209 (2020).
\end{enumerate}

\end{document}